\documentclass[11pt]{amsart}
\usepackage{xypic,amscd,amsmath,amsthm,amssymb,textcomp}
\usepackage{centernot}
\usepackage[dvips]{color}
\newtheorem{theorem}{Theorem}[subsection]
\newtheorem{lemma}[theorem]{Lemma}

\newtheorem{proposition}[theorem]{Proposition}
\newtheorem{corollary}[theorem]{Corollary}

\newcommand{\PP}{\mathbb{P}}

\newcommand{\QQ}{\mathbb Q}

\newcommand{\ZZ}{\mathbb{Z}}

\newcommand{\cE}{\mathcal{E}}

\newcommand{\M}{\mathcal{M}}

\newcommand{\Hom}{{\operatorname{Hom}}}

\newcommand{\n}{\noindent}
\newcommand{\av}{\alpha _{van}}

\newcommand{\frp}{\mbox{$\frac{1}{p}$}}
\theoremstyle{definition}

\theoremstyle{remark}

\begin{document}

\date{}
\title[]
{Some remarks on Brauer Classes of $K3$-type}

\author{Federica Galluzzi}
\address{Dipartimento di Matematica, Universit\`a di Torino, Via Carlo Alberto 10, Torino, Italy}
\email{federica.galluzzi@unito.it}

\author{Bert van Geemen}
\address{Dipartimento di Matematica, Universit\`a di Milano, Via Cesare Saldini 50, Milano, Italy}
\email{lambertus.vangeemen@unimi.it}

%\keywords{Brauer group, Cubic 4-folds, K3 surfaces }

\begin{abstract}
An element in the Brauer group of a general complex projective $K3$ surface $S$ defines a sublattice of the
transcendental lattice of $S$.
We consider those elements of prime order for which this sublattice is Hodge-isometric to the
transcendental lattice of another K3 surface $X$.
We recall that this defines a finite map between moduli spaces of polarized K3 surfaces and we compute its degree.
We show how the Picard lattice of $X$ determines the Picard lattice of $S$
in the case that the Picard number of $X$ is two.
%and
%We obtain some results on possible isomorphisms between these K3 surfaces.
%We also discuss the classical Mukai correspondence between $K3$ surfaces of degree eight and double planes.
\end{abstract}

\maketitle

%\smallskip
%\n 
%{\it Mathematics Subject Classification (2010)}:  14F22, 14E08, 14C30,14J28.

%\tableofcontents

\section*{Introduction}
\setcounter{subsection}{1}
Let $S$ be a complex projective $K3$ surface and let $T(S)$ be its transcendental lattice.
The Brauer group $Br(S)$ of $S$ can be identified with (\cite[\S 18.1]{Hu16})
$$
Br(S)\,=\,\Hom(T(S),\QQ/\ZZ)~.
$$
An element $\alpha$ of order $p$, for a prime number $p$, in the Brauer group $Br(S)$ then defines
a Hodge substructure $T_{\alpha}(S)$ of index $p$ in the transcendental lattice $T(S)$, it is the kernel of $\alpha$.
The isometry classes of these sublattices were determined in \cite{vG05}, \cite{vGK23} for $p=2$
and  in \cite{McK17} for $p>2$. In particular, if the Picard rank of $S$ is one and $Pic(S)=\ZZ h$ with $h^2=2d$,
then for any prime number $p$ there is one class whose lattices are isometric to the transcendental lattice
of a K3 surface of degree $2p^2d$. We will call these {\em Brauer classes of K3-type}.

There is a finite map
$$
\kappa:\;\M_{2p^2d}\,\longrightarrow\,\M_{2d}
$$
between the moduli spaces of polarized K3 surfaces of degree $2p^2d$ and degree $2d$ whose fiber over $S$,
with Picard rank one, consists of all $X$ with $T(X)\cong T_\alpha(S)$ for some $\alpha\in Br(S)_p$ of K3-type.
We determine the degree of this map in Proposition \ref{proposition:Degkappa}.
Since the Hodge structure $T(X)$ does not determine $X$ uniquely in general
due to the presence of Fourier-Mukai partners,
this degree is not simply the number of cyclic subgroups in $Br(S)_p$ of K3-type.

In case the Picard rank of $X$ is two, its Picard lattice is determined by two integers $b,c$ and we write
$X_{b,c}$ for such a surface.
We use the description of $\kappa$
%in terms of moduli spaces of sheaves on $X$ to determine,
to determine the Picard lattice of $S_{b,c}:=\kappa(X_{b,c})$.
In the Picard rank two case it can (and does) happen that $X_{b,c}\cong S_{b,c}$,
if so, the determinants of the Picard lattices are the same. This is an easy necessary,
but not sufficient, criterion for the existence of an isomorphism.
%We refer to \cite{MN03}(??) for the this problem.

The equality of the determinants
occurs exactly when the transcendental lattice of $X_{b,c}$, which is $T_\alpha(S_{b,c})$, is equal to $T(S_{b,c})$.
More precisely, consider a family of K3 surfaces over a disc with special fiber $S_{b,c}$ and general fiber a polarized
K3 surface $S$ of degree $2d$ of Picard rank one. Then we can identify $H^2(S_{b,c},\ZZ)=H^2(S,\ZZ)$
and we have $T(S_{b,c})\subset T(S)$. Therefore there is a restriction map on the Brauer groups and on their
$p$-torsion subgroups. The kernel of this map is the subgroup of vanishing Brauer classes, any non-zero element in
it is called a {\em vanishing Brauer class} (\cite{GvG23}):
$$
\langle \av\rangle\,=\,\ker(Br(S)_p\,\longrightarrow\,Br(S_{b,c})_p)~.
$$
%In \cite{GvG23} we worked this out for the case $p=2$ with applications
%to the K3 surfaces associated to a cubic fourfold.
%we recall these results in \ Section \ref{Section:oddcub}.

Let $X$ be a K3 surface of degree $2p^2d$ with Picard rank one and let $S=\kappa(X)$ so that $T(X)$ has index $p$ in $T(S)$.
We denote by $\alpha_X\in Br(S)_p$ an element such  that $\ker(\alpha_X)=T(X)$.
A specialization of $S$ to $S_{b,c}$ then induces a specialization of $X$ to $X_{b,c}$ and
we have two cyclic subgroups in $Br(S)_p$, one is $\langle\alpha_X\rangle$
and the other is $\langle\av\rangle$. These two subgroups coincide exactly when $\alpha_X$ is trivial on $T(S_{b,c})$,
so when $T_{\alpha_X}(S_{b,c})=T(S_{b,c})$.
%In general, a specialization defines two distinct subgroups of $Br(S)_p$ and make these explicit
In Theorem \ref{classiXS} we make these subgroups explicit.
%in \cite{GvG23} we worked out a classification of the possible pairs of such subgroups.

We intend to use these results to study degree eight K3 surfaces $S_\beta$ associated to certain cubic fourfolds
in \cite[Proposition 5.1.4]{GvG23}. In \cite{KS17}, Kuznetsov and Shinder study the classes generated by
K3 surfaces in the Grothendieck ring of $K_0(\mbox{Var}/{\mathbb K})[{\mathbb L}^{-1}]$.
They use the geometry of the conic bundles associated to a Brauer class of K3-type,
in particular in the case of a specialization in which this class vanishes.
The very basic results in this paper might give some more insight into these cases.

\

\section{Brauer groups, vanishing classes and invariants}\label{brk3}
\setcounter{subsection}{1}

\subsection{Brauer classes and B-fields}
We recall, following \cite{GvG23}, the main definitions but now for the case of an arbitrary prime number $p$
rather than only $p=2$.

Let $S$ be a $K3$ surface, its transcendental lattice is $T(S):=Pic(S)^\perp$ in $H^2(S,\ZZ)$.
The \textit{Brauer group} of $S$ can be identified with  (cf.\ \cite[18.1]{Hu16})
$$
Br (S)= \Hom(T(S),\QQ/\ZZ)~.
$$
Since $H^2(S,\ZZ)$ is a selfdual lattice, any such homomorphism $\alpha$ can be defined by an element $B=B_\alpha\in H^2(X,\QQ)$, called a B-field representative of $\alpha$:
$$
\alpha:\,T(S)\,\longrightarrow\,\QQ/\ZZ,\qquad t\,\longmapsto B\cdot t\mod\ZZ~.
$$

Let $Br(S)_p$ be the $p$-torsion subgroup of $Br(S)$ and let $\alpha\in Br(S)_p$. Then the homomorphism $\alpha$ takes
values in $\frp \ZZ/\ZZ$.
A B-field $B_{\alpha}\in \frp H^2(S,\ZZ)$  is unique
up to $\frp Pic(S) + H^2(S,\ZZ)$ :
$$
B_{\alpha}'=B_{\alpha}+\frp D+c,\qquad D\in Pic(S) ,\quad c \in H^2(S,\ZZ)~,
$$

\subsection{Brauer classes, sublattices and invariants: the case $p=2$}\label{BSIp2}
In case the Picard rank of $S$ is one, the lattices $T_\alpha(S)$ for $S\in Br(S)_p$ are classified up to isometry
by their discriminant groups. This leads to the following classification, for $p=2$ in Lemma \ref{lem:invalpha2}
and for $p>2$ in Lemma \ref{lem:invalphap}.
The case $p=2$ is presented in a format similar to the one for $p>2$.
%Recall that a Brauer class $\alpha$ is of K3-type if $\Gamma_\alpha$

\

\begin{lemma}(\cite[Theorem 2.3]{vGK23}) \label{lem:invalpha2}
Let $S$ be a K3  surface such that $Pic(S)= \ZZ h$, $h^2=2d >0$.
Let $\alpha \in Br(S)_2$, $\alpha\neq 0$, and $B \in \mbox{$\frac{1}{2}$} H^2(S,\ZZ) \subset H^2(S,\QQ)$
a B-field representing $\alpha$.
\begin{enumerate}
\item[a)] In case $2\centernot{|}d$ there are three isomorphism classes of lattices $T_\alpha(S)$.
\begin{enumerate}
\item[i)] $Bh\equiv 0\mod \ZZ$, in this case there is a unique isomorphism class, of order $2^{20}-1$,
with discriminant group $\ZZ/2d\ZZ\oplus\ZZ/2\ZZ\oplus\ZZ/2\ZZ$,
\item[ii)] $Bh\equiv 1/2\mod \ZZ$, $B^2\equiv 0\mod \ZZ$,
in this case there is a unique isomorphism class, of order $2^9(2^{10}+1)$, with discriminant group $\ZZ/8d\ZZ$,
\item[iii)] $Bh\equiv 1/2\mod \ZZ$, $B^2\equiv \mbox{$\frac{1}{2}$} \mod \ZZ$,
in this case there is a unique isomorphism class, of order $2^9(2^{10}-1)$, with discriminant group $\ZZ/8d\ZZ$.
\end{enumerate}
\item[b)] In case $2{|}d$ there are three isomorphism classes of lattices $T_\alpha(S)$.
\begin{enumerate}
\item[i)] $Bh\equiv 0\mod \ZZ$, $B^2\equiv 0\mod\ZZ$, in this case there is a unique isomorphism class of lattices,
of order $2^{9}(2^{10}+1)-1$, with discriminant group $\ZZ/2d\ZZ\oplus\ZZ/2\ZZ\oplus\ZZ/2\ZZ$,
\item[ii)] $Bh\equiv 0\mod \ZZ$, $B^2\equiv \mbox{$\frac{1}{2}$}\mod\ZZ$,
in this case there is a unique isomorphism class of lattices,
of order $2^{9}(2^{10}-1)$, with discriminant group $\ZZ/2d\ZZ\oplus\ZZ/2\ZZ\oplus\ZZ/2\ZZ$,
\item[iii)] $Bh\equiv \mbox{$\frac{1}{2}$}\mod \ZZ$, in this case there is a unique isomorphism class of lattices,
of order $2^{20}$, with discriminant group $\ZZ/8d\ZZ$.
\end{enumerate}
\end{enumerate}

A non-trivial Brauer class $\alpha\in Br(S)_2$ is of K3-type if $B_\alpha h\equiv 1/2\mod\ZZ$ and
$B^2_\alpha\equiv 0\mod \ZZ$ (\cite[Corollary 9.4]{vG05}, the latter is significant only if $2\centernot{|}d$).

In case $d=1$, let $C_6\subset\PP^2$ be the (smooth) degree six branch curve of the double cover $S\rightarrow\PP^2$.
Then the Brauer classes in $Br(S)_2$ correspond to
(\cite{vG05}, \cite{IOOV17}):
\begin{enumerate}
\item[i)] If $B_{\alpha}h \equiv 0$, $\alpha$ corresponds to a point of order two $p \in Jac(C_6)$.
\item[ii)] If $B_{\alpha}h \equiv \mbox{$\frac{1}{2}$}\, $ and $\,B_{\alpha}^2 \equiv 0$, $\alpha$ corresponds to an even theta characteristic on $C_6$.
\item[iii)] If $B_{\alpha}h \equiv \mbox{$\frac{1}{2}$}\,$ and $\,B_{\alpha}^2 \equiv \mbox{$\frac{1}{2}$}$, $\alpha$
corresponds to an odd theta characteristic on $C_6$.
\end{enumerate}
\end{lemma}

\

\subsection{Brauer classes, sublattices and invariants: the case $p>2$}\label{BSIp}
In case the Picard rank of the degree $2d$ K3 surface $S$ is one,
the lattices $T_\alpha(S)$ for $S\in Br(S)_p$ and $p>2$ are classified up to isometry
by their discriminant groups, of order $2p^2d$, $d(T_\alpha(S),q)$,
where $q=q_\alpha$ is a quadratic form with values in $\QQ/2\ZZ$. In case the discriminant group is cyclic,
we denote by $v$ a generator of the discriminant group, $d(T_\alpha(S))=\langle v\rangle$.

To give the classification, we fix an isomorphism $H^2(S,\ZZ)=U^3\oplus E_8(-1)^2$ such that $h=(1,d)\in U$,
the first copy of $U$ in the lattice. Then
$$
T(S)=h^\perp=\ZZ v\oplus\Lambda',\quad  \mbox{with}\; v=(-1,d)\,\in U~.
$$ Let $w:=(0,-1)\in U$. Then  $wv=1$ and
any $\alpha\in \Hom(T(S),\frp\ZZ/\ZZ)$ is determined by a B-field
$B_\alpha=\frp (i_\alpha w+\lambda_\alpha)\in \frp H^2(S,\ZZ)$
where $i_\alpha\in\ZZ$ and $\lambda_\alpha\in \Lambda'$:
$$
\alpha:\,T(S)\,\longrightarrow\, \frp\ZZ/\ZZ,\qquad \alpha(zv+\lambda')\,=\,B_\alpha\cdot(zv+\lambda')\,=\,
\frp (i_\alpha z+ \lambda_\alpha\cdot\lambda')~.
$$
We define $c_\alpha:=-\lambda_\alpha^2/2\in\ZZ$ and we observe that
$$
B_\alpha\cdot h\,=\,-\frp i_\alpha,\qquad \lambda_\alpha^2\,=\,-2c_\alpha~.
$$

\

\begin{lemma} (\cite[Theorem 9]{McK17})\label{lem:invalphap}
Let $S$ be a K3  surface such that $Pic(S)= \ZZ h$, $h^2=2d >0$.
Let $p>2$ be a prime number, $\alpha \in Br(S)_p$ and $B_{\alpha} \in \frac{1}{p} H^2(S,\ZZ) \subset H^2(S,\QQ)$
a B-field representing $\alpha$.

\begin{enumerate}
 \item[a)]
In case $p\centernot{|} d$, there are three isomorphism classes of lattices $T_{\alpha}(S)$,
\begin{enumerate}
\item[i)] the discriminant group is cyclic, so isomorphic to $\ZZ/2p^2d$,
and $-2dp^2q(v)\mod p$ is a quadratic residue, there are $\mbox{$\frac{1}{2}$}p^{10}(p^{10}+1)$ such lattices;
\item[ii)] the discriminant group is cyclic, so isomorphic to $\ZZ/2p^2d$,
and $-2dp^2q(v)\mod p$ is  not a quadratic residue, there are $\mbox{$\frac{1}{2}$}p^{10}(p^{10}-1)$ such lattices;
\item[iii)]  there is a unique isomorphism class of lattices with
discriminant group $\ZZ/2d\ZZ\oplus\ZZ/p\ZZ\oplus\ZZ/p\ZZ$ and there are ${(p^{20}-1)}/{(p-1)}$ such sublattices.
\end{enumerate}
A Brauer class is of K3-type if the discriminant group is cyclic and $-2dp^2q(v)\mod p$ is a square in $\ZZ/p\ZZ$.
\item[b)]
In case $p| d$, there are four isomorphism classes of lattices $T_{\alpha}(S)$,
\begin{enumerate}
\item[i)] $B_\alpha h\equiv 0\mod\ZZ$, the discriminant group is $\ZZ/2d\ZZ\oplus\ZZ/p^2\ZZ$ and $c_\alpha\mod p$ is a
quadratic residue, there are $\mbox{$\frac{1}{2}$}p^{9}(p^{10}-1)$ such sublattices;
\item[ii)] $B_\alpha h\equiv 0\mod\ZZ$, the discriminant group is $\ZZ/2d\ZZ\oplus\ZZ/p^2\ZZ$ and $c_\alpha\mod p$ is not a
quadratic residue, there are $\mbox{$\frac{1}{2}$}p^{9}(p^{10}-1)$ such sublattices;
\item[iii)] $B_\alpha h\equiv 0\mod\ZZ$, there is a unique isomorphism class of lattices with discriminant group is $\ZZ/2d\ZZ\oplus\ZZ/p\ZZ\oplus\ZZ/p\ZZ$, there are $(p^9+1)(p^{10}-1)/(p-1)$ such lattices;
\item[iv)] $B_\alpha h\not \equiv 0\mod\ZZ$, in this case there is a unique isomorphism class of lattices with
discriminant group $\ZZ/2p^2d\ZZ$, there are $p^{20}$ such sublattices.
\end{enumerate}
A Brauer class is of K3-type if the discriminant group is cyclic.
\end{enumerate}
\end{lemma}

\

\subsection{Remarks} The intersection number $B_\alpha h\in \frp \ZZ/\ZZ $ is an invariant of the Brauer class
$\alpha$ only in case $p|d$. In fact, $B_\alpha$ and $B_\alpha+\frp h$ define the same $\alpha$ but
$(B_\alpha+\frp h)h=B_\alpha h +\frp 2d$ which is congruent to $B_\alpha h$ only if $p|d$.

Moreover, if $p|d$ and $B_\alpha h\equiv 0\mod\ZZ$, one obtains the invariant $B_\alpha^2\in (\frp)^2\ZZ/\frp\ZZ$
since any other representative is given by $B_\alpha+\frp D+c$ with $D=ah\in Pic(S)$ and $c\in H^2(S,\ZZ)$.

If $p\centernot{|}d$, we see that there is a choice of the B-field representative such  that
$B_\alpha h=0$, that is, $B_\alpha\in \frp T(S)$, any other such representative is then given by $B_\alpha+c$ with
$c\in T(S)$.

\

\subsection{Vanishing Brauer classes }
Let $S$ be a K3 surface with $Pic(S)=\ZZ h$ and $h^2=2d$. We consider a specialization of $S$ to a K3 surface
$S'$ with Picard lattice
%\begin{equation}\label{pics2}
$$
Pic(S')\,=\,\left(\ZZ h\,\oplus\,\ZZ k,\;
\begin{pmatrix}
h^2 & hk \\
hk & k^2 
\end{pmatrix}   =  \begin{pmatrix}
2d & b \\
b & 2c
\end{pmatrix}
\right)
\quad \mbox{for some \;} b,c \in \ZZ ~.
$$
%\end{equation}

We can identify $H^2(S',\ZZ)=H^2(S,\ZZ)$
and we have $T(S')\subset T(S)$. Thus, there is a restriction map on the
$p$-torsion subgroups of the Brauer groups:
$$
Br(S)_p\,\longrightarrow\,Br(S')_p~.
$$
A non-zero element in the kernel of this map is a {\em vanishing Brauer class} (\cite{GvG23}).

In \cite[Proposition 2.1.2 and Corollary 2.1.3]{GvG23}
we exhibited a B-field representative of a vanishing Brauer class $\av\in Br(S)_2$
which can be easily generalized. To do so we identify $H^2(S,\ZZ)$ with $H^2(S',\ZZ)$
such that $h\in Pic(S)$ specializes to the element with the same name in $H^2(S',\ZZ)$.

\begin{proposition} \label{prop:Brepav}
Let $p$ be a prime number.
We denote by $\av\in Br(S)_p$ a vanishing Brauer class
for the specialization of $(S,h)$ to $S'$ as above.
\begin{itemize}
\item[i)] There is an $\av$ with B-field representative given by
$$
B_{van}\,:=\,\frp k\;(\in \frp H^2(S,\ZZ))~.
$$
\item[ii)] For this $\av \in Br(S)_p$ we have
$$
B_{van} h\,\equiv\, \frp b \mod\ZZ,\quad B_{van}^2\,\equiv\,2c(\frp)^2 \mod \ZZ~.
$$
\end{itemize}
\end{proposition}

\

{\em Proof.}
Notice that $\frp k\not\in \frp Pic(S)+H^2(S,\ZZ)$, hence it defines a non-trivial class in $Br(S)_p$.
But $\frp k\in \frp Pic(S')$ hence $\frp k$ defines the trivial class in $Br(S')_p$.
Therefore the Brauer class with
B-field representative $\frp k$ is the vanishing Brauer class $\av$.
\qed

\setcounter{subsection}{1}

\

\section{K3 surfaces of degree $2d$ and $2p^2d$}\label{Section:K3deg}

\subsection{The Mukai lattice and moduli spaces of sheaves}\label{Mukai}
Let $(X,H)$ be a polarized K3 surface, with $H\in Pic(X)$ primitive,
of degree $H^2=2p^2d$ where $d>0$ and $p$ is a prime number.
As in \cite[\S 2.6, \S 3]{McK17} we consider the Mukai vector
$$
v\,:=\,(p,H,pd)\;\in\,\tilde{H}(X)\,:=\,H^0(X,\ZZ)\oplus H^2(X,\ZZ)\oplus H^4(X,\ZZ)~.
$$
The Mukai lattice $\tilde{H}(X)$ has the bilinear form
$$
(r,c,s)(r',c',s')\,:=\,-(rs'+sr')+c\cdot c',\qquad\mbox{so}\quad v^2\,=\,0~,
$$
where $c\cdot c'$ is the intersection product of $c,c'\in H^2(X,\ZZ)$ and $H^0,H^4$ are naturally identified with $\ZZ$.
The Mukai lattice has the weight two Hodge structure defined by the one on $H^2(S,\ZZ)$. The sublattice
of integral $(1,1)$-classes is thus generated by
the summands $H^0,H^4$ and $Pic(S)\,(\subset H^2$). In particular, $v$ is of type $(1,1)$.

From the work of Mukai \cite{Mu84} it follows that the moduli space $M_X(v)$ of sheaves $\cE$ with
$$
v\,=\,v(\cE)\,:=\,\Big(\mbox{rank}(\cE),c_1(\cE),\mbox{rank}(\cE)+(1/2)c_1(\cE)^2-c_2(\cE)\Big)
$$
is a K3 surface $S$. It is the unique K3 surface for which there is a Hodge isometry
$H^2(S,\ZZ)\cong v^\perp/v$.
This implies that the image of the transcendental lattice $T(X)$ of $X$ under the map
$T(X)\hookrightarrow v^\perp\rightarrow v^\perp/v$ has finite index in $T(S)$.
The Picard ranks of $X$ and $S$ are thus the same.
The sublattice generated by $H^0(X,\ZZ),H^4(X,\ZZ)$ and $\ZZ H$ intersects $v^\perp$ in a rank two sublattice
whose image in $v^\perp/v$ has rank one. Then $(S,h)$, where $h$ is a generator of this rank one lattice,
is a polarized K3 surface of degree $2d$ (cf.\ the proof of Theorem \ref{specPic}).

\subsection{A map between moduli spaces}\label{mapsModuliK3}
This defines a finite map
$$
\kappa\,=\,\kappa_v:\;\M_{2dp^2}\,\longrightarrow\,\M_{2d},\qquad (X,H)\,\longmapsto (M_X(v),h)
$$
where $\M_e$ is the coarse moduli space of K3 surfaces of degree $e$.
The case $d=1$ was used by Kondo in \cite{Ko93}.

\

The geometry behind this map is well understood in the case that $d=1$, $p=2$:
a general K3 surface $(X,H)$ of degree eight determines a K3 surface $(S,h)$ of degree two as follows.
The line bundle $H$ gives an embedding of $X$ as a complete intersection of three quadrics in $\PP^5$.
The surface $S$ is the double cover of the $\PP^2$ that parametrizes the quadrics which is branched over
discriminant curve $C_6\subset\PP^2$ which parametrizes
the singular quadrics (\cite{Kh05}, \cite{IKh13}, \cite{IKh15}, \cite{KS17}, \cite[3.2]{McK17}).

\

%\subsection{The Picard rank two case}

\begin{theorem} \label{specPic}
Let, for a given $d>0$ and prime number $p$, $(X_{b,c},H)$ be a K3 surface of degree $2p^2d$
with Picard lattice
$$
Pic(X_{b,c})\,=\,\left(\ZZ H\oplus\ZZ K,\;
\begin{pmatrix}
2p^2d& b \\
b & 2c
\end{pmatrix} \right ).
$$
Let $S_{b,c}:=\kappa(X_{b,c})$.
The K3 surface $(S_{b,c},h)$ of degree $2d$ has Picard lattice
{\renewcommand{\arraystretch}{1.0}
$$
Pic(S_{b,c})\,=\,\left\{
\begin{array}{rclc}
\left(\ZZ h\oplus\ZZ k,\;
\begin{pmatrix}
2d& b \\
b & 2cp^2
\end{pmatrix} \right )&\mbox{if $p\centernot{|} b$}~,\\
&&&\\
\left(\ZZ h\oplus\ZZ k,\;
\begin{pmatrix}
2d& b/p \\
b/p & 2c
\end{pmatrix} \right )&\mbox{if $p| b$}~.\\
\end{array}\right.
$$
}
\end{theorem}

\noindent {\em Proof.}
First of all we show that the general $S_{b,c}$ has a polarization of degree $2d$.
The sublattice of $(1,1)$ classes in $\tilde{H}(X_{b,c})$ contains the primitive sublattice $N$ generated by
$(1,0,0),(0,H,0),(0,0,1)$. One easily finds that
$$
N\cap v^\perp\,=\,\langle \alpha:=(-1,0,d),\;\beta:=(2p,H,0)\rangle,\qquad v\,=\,p\alpha+\beta~.
$$
Therefore $(N\cap v^\perp)/v\cong \ZZ h$ where $h$ is represented by $\alpha$ and $h^2=\alpha^2=2d$, and $h$
is primitive in $v^\perp/v=H^2(S_{b,c},\ZZ)$, of type $(1,1)$ and $(S_{b,c},h)$ defines a point in $\M_{2d}$.

The Picard lattice of $S_{b,c}$ is the image of $N+\ZZ(0,K,0)$,
the sublattice of all $(1,1)$ classes in $\tilde{H}(X_{b,c})$, in $H^2(S_{b,c},\ZZ)$.
$$
(N+\ZZ(0,K,0))\cap v^\perp\,=\,\langle \alpha,\;\beta,\;\gamma\rangle \quad \mbox{with}\quad
\gamma:=\,\left\{\begin{array}{cr} (0,pK,b),&p\centernot{|}\,b,\\ (0,K,b'),& b=pb'.\end{array}
\right.
$$
As $v=p\alpha+\beta$, the image of this sublattice is generated by the images $h,k$ of $\alpha$ and $\gamma$
respectively and one finds the Gram matrices as in the theorem.
%The subgroup of vanishing Brauer classes has a generator with B-field representative  $\frp k$, see
%Proposition \ref{prop:Brepav}.
\qed

\begin{corollary} \label{corpb}
If the K3 surfaces $X_{b,c}$ and $S_{b,c}$ are isomorphic, then
the prime number $p$ does not divide $b$.
\end{corollary}

\

\noindent{\em Proof.} If the surfaces are isomorphic, the determinants of the Gram matrices of the Picard groups
must be the same. This is the case only if $p$ does not divide $b$.
(In general it is not the case however that if $p$ does not divide $b$ then $X$ and $S$ are isomorphic,
nor that their Picard lattices are isomorphic.)
\qed

%{\it forse dire qualcosa sulla forma `standard' di $Pic(X)$ e $Pic(S)$?}

\

\section{The map $\kappa$ and Brauer groups}

\subsection{The cyclic subgroup $C\subset Br(S)$ determined by $X$} \label{alphaeven}
For a general $(X,H)$, the transcendental lattice $T(X)=H^\perp$  maps to
a sublattice of index $p$ in $h^\perp$, in fact $H^2=p^2h^2$. In particular there is an
isomorphism $T(S)/T(X)\cong\ZZ/p\ZZ$ and hence there is a surjective map $T(S)\rightarrow\frp \ZZ/\ZZ$
whose kernel is $T(X)$. Thus $S=\kappa(X)$ comes with a subgroup $C=C_X\subset Br(S)_p$ of order $p$.

The following theorem identifies the subgroup $C\subset Br(S)_p$.
It also determines the vanishing Brauer class for a specialization
of an $(S,h)$ with Picard rank one to $S_{b,c}$. We recall that one can choose the isomorphism
$H^2(S,\ZZ)\cong U^3\oplus E_8(-1)^2$ in such a way that $h$ maps to the vector $(1,d)$ in the first copy of $U$.
Then $T(S)=h^\perp$ is the sublattice
$$
T(S)\,=\,\ZZ t_S\oplus U^2\oplus E_8^2
\qquad (\subset \Lambda_{K3}\,:=\, U^3\oplus E_8(-1)^2).\qquad
t_S\,:=\,\left (\begin{smallmatrix}-1\\d \end{smallmatrix}\right )~.
$$

In case $X$ has higher Picard rank, let $H\cdot Pic(X)=\gamma\ZZ$. A result of Mukai implies that the index
of $T(X)$ in $T(S)$ is $GCD(p,\gamma)$, as we verify below for the Picard rank two case.
%See Nikulin, Correspondences II for this

\

\begin{theorem} \label{classiXS}
Let $(X,H)$ be a  K3 surface of degree $2p^2d$ with $Pic(X)=\ZZ H$ and let $(S,h)=\kappa(X,H)$.
Then there is an isomorphism $H^2(S,\ZZ)\cong U^3\oplus E_8(-1)^2=U\oplus\Lambda'$ such that
$$
h\,\longmapsto\, \left(\left(\begin{smallmatrix}1\\d \end{smallmatrix}\right ),0\right),\quad
T(S)\,\stackrel{\cong}{\longrightarrow}\,\left(\begin{smallmatrix}-1\\d \end{smallmatrix}\right )\ZZ\oplus \Lambda'~,
$$
and
$$
T(X)\,=\,\ker(\alpha_X:\,T(S)\,\longrightarrow\,\frp\ZZ/\ZZ,\quad t\,\longmapsto B_X\cdot h\mod\ZZ),
$$
where the B-field is
$$
\quad B_X\,:=\,
\frp\left(\left(\begin{smallmatrix}0\\1 \end{smallmatrix}\right ),0\right)\,\in\, (U\oplus\Lambda')\otimes\QQ~.
$$

In the specialization of $S$ to $S_{b,c}$,
%the subgroup of vanishing Brauer classes is generated by $k/2$,
%the second basis vector of $Pic(
the subgroups of $Br(S)_p$ generated by $\alpha_X$ and $\av$ coincide if and only if $b$ is odd.

In case $p=2$, the invariants of $B_X$ are $B_Xh=1/2$ and $B_X^2=0$.
If also $d=1$, $S$ is a double plane, the Brauer class corresponds to an even theta characteristic and
$X$ is a K3 surface of degree eight, moreover,
for the specialization of $S$ to $S_{b,c}$ the vanishing Brauer class $\alpha _{van}$:
\begin{enumerate}
\item[i)]  is $\alpha_X$ and corresponds to an even theta for $b$ odd,
\item[ii)]
corresponds to a theta characteristic for $b\equiv 2\mod 4$ which is even
if $c$ is even, but $\av\not=\alpha_X$, and is odd otherwise,
\item[iii)]
corresponds to a point $p$ of order two in the Jacobian $Jac(C_6)$  for $b\equiv 0\mod 4$.
\end{enumerate}
In case $b\equiv 0\mod 4$
the theta characteristic $\alpha_{van} + \alpha _X$ is even/odd exactly when $c$ is even/odd.
\end{theorem}

\

{\em Proof.}
Up to isometry, there is a unique embedding of $Pic(X_{b,c})$ in the K3-lattice $\Lambda_{K3}=U\oplus\Lambda'$
with $\Lambda'=U^2\oplus E_8(-1)^2$ (\cite[Thm.\ 1.14.4]{Ni80}).
We choose the isometry such that, for some $K'\in \Lambda'$,
$$
H\,=\,((1,p^2d),0)\in U\oplus\Lambda',\qquad K=((0,b),K')\in U\oplus\Lambda',\quad K^2=(K')^2=2c~.
$$
As $Pic(X)=\ZZ H$ we get
$$
T(X)=H^\perp=\ZZ t_X\oplus \Lambda',\quad
t_X\,:=\,\left (\begin{smallmatrix}-1\\p^2d \end{smallmatrix}\right ) ~.
$$
Notice that $\tilde{H}=(H^0(X)\oplus U\oplus H^4(X))\oplus \Lambda'$,
let
$$
\tilde{U}\,:=\,H^0\oplus U\oplus H^4\,\subset\, \tilde{H}(X,\ZZ)~,
$$
where $U$ is the first summand of $\Lambda_{K3}$. Then $v=(p,H,pd)\in\tilde{U}$.
As $\Lambda'\subset v^\perp$ and $\langle v\rangle \cap \Lambda'=\{0\}$,
this unimodular lattice maps isomorphically to the sublattice
$\Lambda'\subset v^\perp/v=H^2(S,\ZZ)$. To find the image of $T(X)$,
it remains to find the image in $v^\perp/v$ of $(0,(-1,p^2d),0)\in \tilde{U}\cap v^\perp$.

With the notation in the proof of Theorem \ref{specPic},
$$
\tilde{U}\cap v^\perp\,=\,\langle \alpha=(1,0,-d),\;
\beta_1=(0,\left(\begin{smallmatrix}1\\0 \end{smallmatrix}\right ),pd),\;
\beta_2=(0,\left(\begin{smallmatrix}0\\p \end{smallmatrix}\right ),1) \rangle~.
$$
Notice that
$$
v\,=\,p\alpha\,+\,\beta_1\,+\,pd\beta_2,\qquad (0,(-1,p^2d),0)\,=\,-\beta_1+pd\beta_2~.
$$
Hence $(\tilde{U}\cap v^\perp)/v$ is generated by the images $h,k'$ of $\alpha$ and $\beta_2$ whereas
$\beta_1$ maps to $-ph-pdk'=-p(h+dk')$. The intersection products are $h^2=2d$, $hk'=-1$, $(k')^2=0$.
The sublattice $\langle h,k'\rangle\subset H^2(S,\ZZ)$ is isomorphic to $U$:
$$
\langle h,k'\rangle\,\stackrel{\simeq}{\longrightarrow} \,U,\qquad
h\,\longmapsto (1,d),\quad k'\,\longmapsto (0,-1)~.
$$
Then $(0,(-1,p^2d),0)=-\beta_1+pd\beta_2$ maps to $p(h+dk')+pdk'=p(h+2dk')$ which maps to $p(1,-d)\in U$.
So we can choose the isomorphism between $H^2(S,\ZZ)$ and $\Lambda_{K3}$ in such a way that $h\mapsto (1,d)$
and the image of $(0,(-1,p^2d),0)$ maps to $p(1,-d)$.
The image of $T(X)$ in $H^2(S,\ZZ)$ is then the sublattice
$$
T(X)\,\stackrel{\simeq}{\longrightarrow}\,
p\left(\begin{smallmatrix}-1\\d \end{smallmatrix}\right )\ZZ\,\oplus \,\Lambda'\quad(\subset T(S))~.
$$
Notice that for $B_X=(0,1)\in U$, we have $B_X(1,-d)=-1$ which implies that $\ker(B_X)=T(X)$.

A simple computation shows that $B_Xh=\frp (0,1)\cdot(1,d)=\frp$ and $B_X^2=0$.

The subgroups generated by $\alpha_X$ and $\av$ are the same if and only if
the inclusion $T({X_{b,c}})\subset T({S_{b,c}})$ is an equality,
which is equivalent to these lattices having the same discriminants. This is again equivalent to
the determinants of the Picard lattices of $X_{b,c}$ and $S_{b,c}$ being the same and by
Theorem \ref{specPic} we see that this happens if and only if $p\not |\,b$.

The invariants of $\av$ are determined by the second column of a Gram matrix of $Pic(S_{b,c})$ by Proposition \ref{prop:Brepav}. A Gram matrix is given in Theorem \ref{specPic} and (i)-(iii) follow.

In case $B_{van}$ corresponds to a point of order two, the sum $B_s:=B_X+B_{van}$
corresponds to a theta characteristic. The parity of this characteristic is determined by $B_s^2\mod \ZZ$.
We find $B_X^2=0$ and,
by Proposition \ref{prop:Brepav}, $B_{van}=(1/2)k$ with $k$ as in Theorem \ref{specPic},
hence $B_{van}^2=(1/4)(2c)=c/2$. It remains to compute $2B_X\cdot B_{van}$ which we claim is $0$, so that
$B_s^2=c/2$ and the last statement of the theorem is proven.

To verify the claim, we recall that $k$ has representative $\gamma\in v^\perp$ and since $b\equiv 0\mod 4$ we have
$\gamma=(0,K,b')$ where $b=2b'$. With our choice of embedding, $K=((0,b),K')$ and then
$\gamma=(0,(0,b),b')+(0,K',0)=b'\beta_2+(0,K',0)$ which maps to $b'(0,-1)+K'\in U\oplus\Lambda'=H^2(S,\ZZ)$.
Hence $B_X\cdot B_{van}\,=\,(0,1)\cdot (0,-b'/2)=0$ since $B_X\cdot K'=0$.
\qed

\

\subsection{An intrinsic description of the map $\kappa$}\label{kappaintrin}
To define $\kappa:\M_{2p^2d}\rightarrow \M_{2d}$, we used a Mukai vector $v$. Here we give
another way to define the map, where we use some of the notation from the (proofs of the) previous results.

Let $(X,H)$ be a polarized K3 surface of degree $2p^2d$.
Recall that in $H^2(X,\ZZ)$ we have the sublattices $\ZZ H$ and $H^\perp=\ZZ t_X\oplus \Lambda'$, their direct
sum has index $2p^2d$ in $H^2(X,\ZZ)$. To get all of the second cohomology group one has to add the `glue vector'
$(H+t_X)/2p^2d$. Since the discriminant group $(H^\perp)^*/H^\perp$ of $H^\perp$ is cyclic,
there is a unique overlattice, denoted by $h^\perp$, such that $H^\perp\subset h^\perp$ has index $p$.
This overlattice is generated by $H^\perp$ and a $t_S\in h^\perp$ with $pt_S=t_X$.
Let $\ZZ h$ be the rank one lattice with $h^2=2d$. Then the overlattice of $\ZZ h\oplus h^\perp$ defined by the
glue vector $(h+t_S)/2d$ is an even unimodular lattice and hence is isometric to $\Lambda_{K3}$.
This lattice has the Hodge structure induced by the one on $T(X)\subset h^\perp$ and hence defines a unique
polarized K3 surface $(S,h):=\kappa((X,H))$ by surjectivity of the period map and the Torelli theorem.

\

\subsection{FM partners}\label{FMpartners}
We consider the fibers of the map $\kappa$ in the case $d=1,p=2$.
Given a general K3 surface $(S,h)\in \M_2$, by Theorem \ref{classiXS} an $(X,H)$ in the fiber over it
determines an order $2$ subgroup $C=\ker(\alpha_X) \subset Br(S)_2$.
Moreover, the unique non-trivial Brauer class $\alpha\in C$ corresponds to an even theta characteristic on $C_6$.
Given an even theta characteristic on a general $C_6$, this invertible sheaf has no non-trivial global sections and
using \cite{Be00} one obtains a K3 surface $X$ of degree $8$ from a resolution of this sheaf.
Since there are $2^{9}(2^{10}+1)$ even theta characteristics on the genus $10$ curve
$C_6$, this number is also the degree of $\kappa:\M_8\rightarrow\M_2$.

\

For a general $d\geq 1$ and a prime number $p$ however,
the order $p$ subgroup $C=\ker(\alpha_X)\subset Br(S)_p$ only determines the sublattice
$$
T_C\,=\,T_{\alpha_X}(S)\,:=\,\ker(\alpha_X :\,T(S)\,\longrightarrow\,\frp\ZZ/\ZZ)~,
$$
with the induced Hodge structure from $T(S)$.
To obtain a K3 surface $X$, one must embed $T_C$ primitively into a K3 lattice,
which can be done only if $\alpha_X$, and thus $C$, is of K3-type,
and even then the embedding need not be unique up to isometry.

\begin{proposition}\label{proposition:Degkappa}
Let $(S,h)\in\M_{2d}$ be a polarized K3 surface with $Pic(S)=\ZZ h$ and $h^2=2d$. For a prime number $p$,
the cardinality of the fiber of $\kappa:\M_{2p^2d}\rightarrow\M_{2d}$ over $(S,h)$ is
{\renewcommand{\arraystretch}{1.5}
$$
\sharp \kappa^{-1}(S,h)\,=\,\left\{\begin{array}{cl}
                                    \mbox{$\frac{1}{2}$}p^{10}(p^{10}+1) &\mbox{if}\quad d=1,\\
                                    p^{10}(p^{10}+1)&\mbox{if}\quad p\centernot{|}d,d>1,\\
                                    p^{20}&\mbox{if}\quad p|d~.
                                   \end{array} \right.
$$
}
A polarized K3 surface $(X,H)$ is in the fiber $\kappa^{-1}(S,h)$ if and only if $T(X)$ is Hodge isometric to a
sublattice of index $p$ of $T(S)$.
\end{proposition}

\n
{\em Proof.}
If $(X,H)\in \kappa^{-1}(S,h)$ then $T(X)$ is an index $p$ sublattice of $T(S)$,
hence it is defined by an order $p$ subgroup of K3-type of $Br(S)_p$.
Conversely such a subgroup $C$ defines a sublattice $T_C$
of index $p$ that can be primitively embedded into the K3 lattice,
that is, this sublattice is isometric to $H^\perp$ for some $(X,H)$ if $C$ is of K3-type.

The `forgetful' map that associates to a polarized K3 surface $(X,H)$, with $H^2=2p^2d$,
the Hodge structure $H^\perp$ defines a map $\M_{2p^2d}\rightarrow \M^T_{2p^2d}$
which is finite and has degree equal to the number of FM partners of $X$ for a K3 surface $X$ with $Pic(X)=\ZZ H$.
This number is (\cite[Proposition 1.10]{Og02})
$$
|FM(X)|\,=\,2^{\tau(p^2d)-1},
$$
where $\tau(p^2d)$ is the number of prime factors of $p^2d$.
Thus, given the Hodge structure $T_C$, there are $2^{\tau(p^2d)-1}$ K3 surfaces $X$ with $T(X)\cong T_C$.
However, there are also $2^{\tau(d)-1}$ K3 surfaces $S'$ with $T(S')\cong T(S)$.

In particular, if $p|d$ then $\tau(d)=\tau(p^2d)$ and thus, given $T_C$, each $S'$ with $T(S')\cong T(S)$
determines a unique $X'$ with $T_C\cong T(X')$.
The same is true if $d=1$: $\tau(1)=\tau(p^2)=1$.

If however $p\centernot{|}d$ and $d>1$ then $\tau(p^2d)=\tau(d)+1$ and $S'$ determines two K3 surfaces $X'$.

The degree of the map, for $p=2$, $p>2$ then follows from \cite[Theorem 2.3]{vGK23} and \cite[\S 2.6]{McK17}
respectively, which gives the number of $C$ of K3-type.
(In \cite[\S 2.7, Remark]{McK17},
it is stated that, based on a result of Kondo, for $d=1$ the degree of $\kappa$ is $p^{10}(p^{10}+1)$,
however Kondo \cite{Ko24} confirmed that the degree is $\mbox{$\frac{1}{2}$}p^{10}(p^{10}+1)$.)
\qed

\

\

\section{Examples}

\subsection{K3 surfaces with a line}
We consider some special cases of $\kappa:X_{2p^2d}\mapsto S_{2d}$
with Picard rank two. In the literature we found the cases
$$
(d,p)\,=\,(1,2),\quad (1,3),\quad (2,2),\quad (2,3),\quad (3,2)~.
$$
We consider these cases where moreover $X$ has a line.
Then $X=X_{1,-1}$ in the notation of Theorem \ref{specPic}, so $b=1$.
The Picard lattices of $X_{1,-1}$ and $S_{1,-1}$ have the same determinant by Corollary \ref{corpb}.
We then consider whether the K3 surfaces $X_{1,-1}$ and $S_{1,-1}$ are isomorphic,
which is not always the case as we will see.

\n
Since the determinants are the same, we have the equality of the Brauer classes $\alpha_X=\av\in Br(S)$,
that is, the Bauer class $\alpha_X\in Br(S_{1,-1})$ defining $X_{1,-1}$ (up to possible FM partners) is trivial.
It would be interesting to see explicitly a rational section of the specialization of the
conic bundles with Brauer class $\alpha_X$ on a general $S$ to $S_{1,-1}$
in this case. Such conic bundles are quite explicitly known in some of the examples.

\

\subsection{$X_8$ and $S_2$ : square determinants}\label{sqdet}
We consider first the K3 surfaces $X_{b,c}$ of degree $8=2p^2d$ with $d=1,p=2$,
such that the Picard lattice has determinant $D$ which is a square: $D=-\det(Pic(X_8))=b^2-16c$.
In that case $X_{b,c}$ has a genus one fibration given by a divisor class $E$ with $E^2=0$.

In \cite[Lemma 3.10]{KS17} it is shown that there are infinitely many $D=-\det(Pic(X_8))=b^2-16c$ for which
$X_{b,c}$ and $S_{b,c}$, with $b$ odd, are not isomorphic, even if the Picard lattices have the same determinant.
For this they consider the case  $D=m^2$ for an odd integer $m$,
for example take $b=m$ and $c=0$.
There is an isomorphism $X_{b,c}\cong S_{b,c}$ if and only if $r^2-Ds^2=\pm 8$ has an integer solution according
to \cite{MN03}. This can now be written as $r^2-(ms)^2=\pm 8$ and if one takes $ms>3$ then $|r^2-(ms)^2|>8$ for any
$r\neq \pm ms$, hence the result.

A simple geometric example where these surfaces are not isomorphic is thus the case that $X$ is a smooth complete intersection of three quadrics that contains a rational normal cubic curve.
Then $X=X_{3,-1}$ and $K$ is the class of the cubic curve:
$$
Pic(X_{3,-1})\,=\,\left(\ZZ H\oplus\ZZ K,\;
\begin{pmatrix}
8& 3 \\
3 & -2
\end{pmatrix} \right ),\qquad D=-\det(Pic(X_8))=25,
$$
and $X_{3,-1} \ncong S_{3,-1}$, in fact, the Picard lattices are not isometric. From Theorem \ref{specPic}
we have
$$
Pic(S_{3,-1})\,=\,\left(\ZZ h\oplus\ZZ k,\;
\begin{pmatrix}
2& 3 \\
3 & -8
\end{pmatrix} \right )\,\cong\,\left(\ZZ h\oplus\ZZ e,\;
\begin{pmatrix}
2& 5 \\
5 & 0
\end{pmatrix} \right ),\quad e:=h+k~.
$$
If the Picard lattices were isometric, there should also be a $(-2)$-class, like $K$, in $Pic(S_{3,-1})$.
However, $(xh+ye)^2=2x^2+10xy=2x(x+5y)$ and this is $-2$ only if either $x=1$, $x+5y=-1$ or $x=-1$, $x+5y=1$,
however both are impossible for $(x,y)\in \ZZ^2$.

On the other hand, in the cases
$D=1,9$, one can take $X_{1,0}$ and $X_{3,0}$ respectively and these are isomorphic to $S_{1,0}$, $S_{3,0}$
respectively since $(r,s)=(\pm 3, \pm 1)$, $(r,s)=(\pm 1, \pm 3)$ give solutions to $r^2-Ds^2=\pm 8$.
These two cases are well-known.

\n
For $D=1$ the two Picard lattices are both isomorphic to
the hyperbolic plane $U$ and $X_{1,0} \cong S_{1,0}$ since the glueing of $U$ to $T(X)=T(S)$ is just a direct sum.
These surfaces have a unique elliptic fibration.
(see \cite[Proposition 3.2.1]{MN03}).

\n
In the case $D=9$ the surface $X_{3,0} \cong S_{3,0}$ is a $K3$ surface of bidegree $(2,3)$ in
$\mathbb P^1 \times \mathbb P^2$.
See \cite[Proposition 3.2.1]{MN03}, \cite[\textsection 5.8]{vG05},
\cite{Be22}),\cite[Proposition 3.7]{IKh13}).

The classical association $X_8\mapsto S_2$, already mentioned in \S \ref{FMpartners},
is studied in (the list is surely not complete)
 \cite{IKh13}, \cite{IKh15}, \cite{Kh05}, \cite{KS17}, \cite{MN03}, \cite{McK17}.

\

\subsection{$X_{16}$ with a line and $S_4$ }
For $d=2$, $p=2$ one has $X_{2p^2d}=X_{16}$. We assume that this surface contains a line $L$.
Then $X_{16}=X_{1,-1}$ and with $K$ the class of the line:
$$
Pic(X_{1,-1})\,=\,\left(\ZZ H\oplus\ZZ K,\;
M_{16}\,=\,\begin{pmatrix}
16 & 1 \\
1 & -2
\end{pmatrix} \right ).
$$

\n
From Theorem \ref{specPic} one finds
\n
$$
Pic(S_{1,-1})\,=\,\left( \ZZ h\,\oplus\,\ZZ k,\;\begin{pmatrix}
                      4&1\\1&-8
                     \end{pmatrix}\right).
$$

\n
The two Picard lattices are isomorphic:
$$
\left\{ \begin{array}{rcl} h&=&H-2K,\\ k&=&-H+3K,\end{array}\right.
\qquad
\left\{ \begin{array}{rcl} H&=&3h+2k,\\ K&=&h+k,\end{array}\right.
$$
and $\det(Pic(X_{16}))=\det(Pic(S_4))=-33$.
To show that $X_{1,-1}\cong S_{1,-1}$ it suffices to show that the glueing of the Picard lattice to the
transcendental lattice is unique up to isomorphisms. That again follows from the surjectivity of the map
$O(Pic(X_{1,-1}))\rightarrow O(D_P)$ where $D_P$ is the discriminant group of $Pic(X_{1,-1})$.
In fact, let
$$
S\,:=\,\begin{pmatrix}
                      19&64\\8&27
                     \end{pmatrix}, \qquad \mbox{then}\quad
                     SM_{16}S^t\,=\,M_{16}~,
$$
so that $S\in O(Pic(X_{16}))$. The discriminant group of the Picard group is generated by $\delta:=(2,1)/33$
%and $q(\delta=2/33$.
and one finds that $\delta S\equiv 23\delta$. Since %$(23)^2\equiv 1\mod 33$ and
$\ZZ/33\ZZ\cong \ZZ/3\ZZ\times \ZZ/11\ZZ$, with $23\mapsto (-1,1)$,
we see that $-I,S\in O(Pic(X_{16}))$ generate $O(D_P)$.

Geometrically, the isomorphism $X_{1,-1}\rightarrow S_{1,-1}$ is given by the `double projection' from the line
$L\subset X_{1,-1}\subset \PP^{9}$. First one projects from the line: $\phi_{H-K}:X_{16}\rightarrow X'_{12}\subset \PP^7$, notice that $(H-K)^2=12$. The image of $L$ is a rational normal curve of degree $(H-K)K=3$
which spans a $\PP^3\subset\PP^7$.
Projection from the span of the normal curve induces the map $\phi_{H-2K}:X_{16}\rightarrow S_4\subset\PP^3$,
the image of $L$ is a quintic rational curve in the quartic surface $S_4$ since $(H-2K)K=5$.

See \cite{IR05}, \cite{IR07}, \cite[\S 3.4]{McK17}, \cite[\S 5.3]{vGK23} for geometrical aspects of the map
$X_{16}\mapsto S_{4}$.

\subsection{$X_{18}$ with a line and $S_2$}   For $d=1$, $p=3$ one has $X_{2p^2d}=X_{18}$ and $S_{2d}=S_2$.
Assume that $X$ contains a line, then $X=X_{1,-1}$ and
$$
Pic(X_{1,-1})\,=\,\left(\ZZ H\oplus\ZZ K,\;
\begin{pmatrix}
18 & 1 \\
1 & -2
\end{pmatrix} \right )~,
$$
therefore
$$
Pic(S_{1,-1})\,=\,\left( \ZZ h\,\oplus\,\ZZ k,\;\begin{pmatrix}
                      2&1 \\1 \, &-18
                     \end{pmatrix}\right).
$$
The two Picard lattices are isomorphic and have determinant $-37$,
%$\det(Pic(X_{18})=\det(Pic(S_2))=-37$,
for example an isomorphism is:
$$
\left\{ \begin{array}{rcl} h&=&2H-5K,\\ k&=&-5H+13K,\end{array}\right.
\qquad
\left\{ \begin{array}{rcl} H&=&13h+5k,\\ K&=&5h+2k~.\end{array}\right.
$$
The two $K3$'s are also isomorphic since the  (sufficient) conditions in \cite[Theorem 3.1.5]{MN04} are satisfied.
More precisely, there exists $h_1 \in Pic(X_{18})$ such that
$h_1^2 =  2 p $  and $h_1 H \equiv 0 \mod p$  (here $p=3$), for example $h_1=H+3K$.

Another way to see this is to notice that, since the order of the discriminant groups is $37,$ a prime number,
the orthogonal group of the discriminant lattice of the Picard groups is $\{\pm 1\}$.
Thus the glueing of the Picard lattice to the transcendental lattice is unique and the surfaces are isomorphic.

See also \cite[\S 3.3]{McK17} for the map $X_{18}\mapsto S_2$.

\subsection{$X_{24}$ with a line and $S_6$}
For $d=3$, $p=2$ one has $X_{2p^2d}=X_{24}$ and $S_{2d}=S_6$.
Assume that $X$ contains a line, then $X=X_{1,-1}$ and the Picard lattices of $X$ and $S_{1,-1}$ have the Gram matrices
$$
P_{24}\,=\,
\begin{pmatrix}
24 & 1 \\
1 & -2
\end{pmatrix},\qquad
P_6\,=\,\begin{pmatrix}
                      6&1 \\1 \, &-8
                     \end{pmatrix}.
$$
The determinants are $-49$, so we are in the case of square determinants as in \S \ref{sqdet}.
Let $e=h+k$ in $Pic(S_{1,-1})$, then $e^2=(h+k)^2=6+2-8=0$. Then $e=h+k, h$ is a basis of the Picard lattice
of $S_{1,-1}$
and
$$
(xe+yh)^2\,=\,x^2e^2+2xyeh+y^2h^2\,=\,14xy\,+\,6y^2\qquad (e=h+k)~.
$$
Notice that there is no $(-2)$-vector in this lattice since $y(7x+3y)=-1$ has no integer solutions.
Therefore the Picard lattices are not isometric and hence $X_{1,-1}\not\cong S_{1,1}$.

The case $X_{24}\mapsto S_6$ was studied in detail in the recent paper \cite{KM23}.

\subsection{$X_{36}$ with a line and $S_4$} For $d=2$ and $p=3$ one has $X_{2p^2d}=X_{36}$, which has genus $19$,
and $S_{2d}=S_4$.
This case was considered in \cite[\S 1.2]{KM23}, see also \cite[Remark 4.16]{BBFM23}.
Assume that $X$ contains a line, then $X=X_{1,-1}$ and the Picard lattices of $X$ and $S_{1,-1}$ have the Gram matrices
$$
P_{36}\,=\,
\begin{pmatrix}
36 & 1 \\
1 & -2
\end{pmatrix},\qquad
P_4\,=\,\begin{pmatrix}
                      4&1 \\1 \, &-18
                     \end{pmatrix}.
$$
The determinants are $-73$. The Gram matrices are equivalent:
$$
SP_{36}S^t\,=\,P_4,\qquad S\,=\,\begin{pmatrix}
57 & 272 \\
136 & 649
\end{pmatrix}~,
$$
hence the Picard lattices are isomorphic. Since the orthogonal group of the discriminant group, which is $\ZZ/73\ZZ$,
is $\{\pm 1\}$, the glueing of the Picard lattice to the transcendental lattice is unique.
%there is a unique K3 surface with this Picard lattice and a given transcendental lattice.
Thus $X_{1,-1}\cong S_{1,-1}$.

\bibliographystyle{amsplain}

\end{document}